\documentclass[12pt]{amsart}
\usepackage{latexsym, amssymb, amsmath}
 \usepackage{eucal}

\title[Order Polynomials and Their Generalizations]
{A New Approach to Order Polynomials of Labeled Posets and Their Generalizations}
    \author{John Shareshian, David Wright and Wenhua Zhao*}      
    
\begin{document}
\keywords{Labeled posets, order polynomials,  
Eulerian polynomials, Bernoulli numbers and Bernoulli polynomials}
   
\subjclass[2000]{06A07, 06A11, 11B68}

    \bibliographystyle{alpha}
    \maketitle

    \newtheorem{rema}{Remark}[section]
    \newtheorem{propo}[rema]{Proposition}
   \newtheorem{theo}[rema]{Theorem}
 \newtheorem{conj}[rema]{Conjecture}
\newtheorem{Algo}[rema]{Algorithm}
   \newtheorem{defi}[rema]{Definition}
    \newtheorem{lemma}[rema]{Lemma}
    \newtheorem{corol}[rema]{Corollary}
     \newtheorem{exam}[rema]{Example}
 \newtheorem{theo-def}[rema]{Theorem-Definition}
  \newtheorem{rmk}[rema]{Remark}
\newtheorem{quest}[rema]{Question}
      \newtheorem{redu}[rema]{Reduction Conditions}
	\newcommand{\nno}{\nonumber}
	\newcommand{\lbar}{\big\vert}
        \newcommand{\mbar}{\mbox{\large $\vert$}}
	\newcommand{\p}{\partial}
	\newcommand{\dps}{\displaystyle}
	\newcommand{\bra}{\langle}
	\newcommand{\ket}{\rangle}
\newcommand{\kr}{\mbox{\rm Ker}\ }
 \newcommand{\res}{\mbox{\rm Res}}
\renewcommand{\hom}{\mbox{\rm Hom}}
 \newcommand{\pf}{{\it Proof:}\hspace{2ex}}
 \newcommand{\epf}{\hspace{2em}$\Box$}
 \newcommand{\epfv}{\hspace{1em}$\Box$\vspace{1em}}
\newcommand{\nord}{\mbox{\scriptsize ${\circ\atop\circ}$}}
\newcommand{\wt}{\mbox{\rm wt}\ }
\newcommand{\clr}{\mbox{\rm clr}\ }
\newcommand{\ideg}{\mbox{\rm Ideg}\ }
\newcommand{\GC}{G_{\mathbb C}}
\newcommand{\gC}{{\mathfrak g}_{\mathbb C}}
\newcommand{\hatC}{\widehat {\mathbb C}}
\newcommand{\bC}{{\mathbb C}}
\newcommand{\bZ}{{\mathbb Z}}
\newcommand{\bQ}{{\mathbb Q}}
\newcommand{\bR}{{\mathbb R}}
\newcommand{\bN}{{\mathbb N}}
\newcommand{\bT}{{\mathbb T}}
\newcommand{\cI}{{\mathcal I}}
\newcommand{\cIw}{{\mathcal I}_\omega}
\newcommand{\cG}{{\mathcal G}}
\newcommand{\Pw}{ {P, \omega} }
\newcommand{\BQ}{\begin{eqnarray}}
\newcommand{\EQ}{\end{eqnarray}}
\newcommand{\BQn}{\begin{eqnarray*}}
\newcommand{\EQn}{\end{eqnarray*}}
\newcommand{\wtilde}{\widetilde}
\newcommand{\Hol}{\mbox{Hol}}
\newcommand{\Hom}{\mbox{Hom}}

\begin{abstract}
In this paper, we first give formulas
for the order polynomial $\Omega (\Pw; t)$ 
and the Eulerian polynomial 
$e(\Pw; \lambda)$ of a finite 
labeled poset $(P, \omega)$
using the adjacency matrix
of what we call the $\omega$-graph of $(P, \omega)$. 
We then derive various recursion formulas for
$\Omega (\Pw; t)$ and $e(\Pw; \lambda)$ and 
discuss some applications of these formulas
to Bernoulli numbers
and Bernoulli polynomials.
Finally, we give a recursive algorithm using a single linear
operator on a vector space. This algorithm
provides a uniform method to construct 
a family of new invariants for labeled posets $(\Pw)$, which 
includes the order polynomial $\Omega (\Pw; t)$ and the
invariant $\tilde e(\Pw; \lambda)
=\frac {e(\Pw; \lambda)}{(1-\lambda)^{|P|+1}}$.
The well-known quasi-symmetric function 
invariant of labeled posets 
and a further generalization of our construction 
are also discussed.
\end{abstract}

\tableofcontents

\renewcommand{\theequation}{\thesection.\arabic{equation}}
\renewcommand{\therema}{\thesection.\arabic{rema}}
\setcounter{equation}{0}
\setcounter{rema}{0}
\setcounter{section}{0}

\section{\bf Introduction}\label{S1}

In this paper, we study order polynomials $\Omega (\Pw; t)$ and Eulerian
polynomials $e(\Pw; \lambda)$ of labeled posets $(\Pw)$
(see \cite{St1} and \cite{St2}) and their generalizations. 
Motivated by the way that strict order polynomials
$\bar \Omega (T)$ of rooted trees appear as coefficients
in the tree expansions of formal flows of analytic maps of $\bC^n$ 
in \cite{WZ}, 
we first give formulas in Section \ref{S2} for  
the order polynomial $\Omega (\Pw; t)$ and the Eulerian
polynomial $e(\Pw; \lambda)$ in terms of 
the adjacency matrix $A (\Pw)$ of what we call the 
$\omega$-graph of a labeled poset $(\Pw)$. 
We also derive various recursion formulas 
for $\Omega (\Pw; t)$ and $e(\Pw; \lambda)$. 
Some applications of these formulas to Bernoulli numbers 
and Bernoulli polynomials are discussed in Section \ref{S2.3}. 
Although some results in Section \ref{S2} can be found in literature,
it seems to us that there is no source that treats 
order polynomials and Eulerian polynomials of labeled posets 
in the manner presented here.
In Section \ref{S3}, motivated by the recursion formulas 
derived in Section \ref{S2} and also the construction of
a family of invariants of rooted trees in \cite{Z}, 
we construct a family of 
invariants for labeled posets
by a recursive algorithm (see Algorithm \ref{Algo}) in terms of
a single linear operator $\Xi$ on a vector space $A$. 
We then show that the order polynomial $\Omega (\Pw; t)$ and the
invariant $\tilde e(\Pw; \lambda)
=\frac {e(\Pw; \lambda)}{(1-\lambda)^{|P|+1}}$ of a labeled poset 
$(\Pw)$ belong to this family of invariants. 
In Section \ref{S3.3},  we consider the 
quasi-symmetric function 
invariant $K(\Pw; x)$ (see \cite{G}, \cite{MR} and \cite{St3})
of a labeled poset $(\Pw)$. Even though  
$K(\Pw; x)$ can not be recovered from our general construction, 
we give a recursion formula for $K(\Pw; x)$ 
(see Proposition \ref{Recur-K}) which suggests 
a further generalization of our construction. Finally, we consider 
invariants of  unlabeled posets $P$ derived  
from Algorithm \ref{Algo} by identifying unlabeled posets 
with labeled posets in certain ways.

\renewcommand{\theequation}{\thesection.\arabic{equation}}
\renewcommand{\therema}{\thesection.\arabic{rema}}
\setcounter{equation}{0}
\setcounter{rema}{0}

\section{\bf Order Polynomials of Labeled Posets}\label{S2}

Once and for all, we fix the following notation and  conventions. 
\vskip2mm

{\bf Notation and Conventions:}
\begin{enumerate}

\item We denote by $\bN$ (resp. $\bN^+$) the set of
all non-negative (resp. positive) integers. For any $n\in \bN^+$, 
the totally
ordered poset $\{1, 2, \cdots, n\}$ will be denoted by $[n]$.
 
\item For any $n\in \bN$, we denote by $\mathbb P_n$ the set of posets with $n$ elements.
We set $\mathbb P=\bigsqcup_{n\in \bN} \mathbb P_n$. 
The empty set $\emptyset$ will always be treated as a finite poset 
and will be viewed as an ideal of 
any poset $P$. The poset with a single element is called the {\it singleton} 
and is denoted 
by $\bullet$. 

\item For any poset $P\in \mathbb P$, we denote by $V(P)$ the set of 
elements of $P$, and by
$|P|$ the cardinality of $V(P)$. 
We denote by $L(P)$ the set of minimum elements of $P$.

\item For any $P\in \mathbb P$, we denote by $\mathcal I(P)$ the set of
all ideals including $\emptyset$ and $P$ itself.

\item For any $n\in \bN$, we denoted by $A_n$ 
the anti-chain with $n$ elements.
The empty poset  $\emptyset$
and the singleton $\bullet$ are anti-chains.

\item We say  a poset $P$ is a {\it
rooted tree} if its Hasse diagram is a rooted tree 
with the root being 
the unique minimum element.

\item  For any $n\in \bN$,  we let
$S_n$ be the {\it shrub} with $n+1$ vertices. 
It is the poset in $\mathbb P_{n+1}$ whose
Hasse diagram is the
unique rooted tree of height $\le 1$. (Equality holds
unless $n=0$.)  
Note that $S_0=\bullet$, the singleton.


\item For a finite poset $P$, we set the following notation: 
\begin{enumerate}
\item $S\trianglelefteq P$: $S$ is an ideal of $P$. 
\item $S\vartriangleleft P$: $S\trianglelefteq P$ but  $S\neq P$.
\end{enumerate}
\end{enumerate}

In this paper ``poset" will always mean a partially ordered finite set.
Recall that a labeled poset $(P, \omega)$ is a poset $P$ with an injective map
$\omega: V(P)\to \bN^+$. When $\omega$ as a map from $P$ to $\bN^+$ is 
order-preserving (resp. order-reversing), we say $(\Pw)$ is a
{\it naturally} (resp. {\it strictly}) {\it labeled} poset. 
We say two labeled posets $(P, \omega)$ and $(Q, \phi)$ 
are {\it equivalent} 
if there is an isomorphism $f: P\to Q$ of posets and an order-preserving
bijection
$g: \text{Im}(\omega)\to \text{Im}(\phi)$ such that 
$g\circ \omega=\phi\circ f$. Note that
equivalent posets have same order polynomial. 
In this paper, the term ``labeled poset" will always refer to an
equivalence class.
In other words, any two equivalent labeled posets 
are considered to be the same poset. Note that each 
equivalence class has a 
unique representative $(\Pw)$ with $\omega(P)=[|P|]$.

Given $(\Pw)$ a labeled poset and $S$ a subset of $V(P)$,
$(S, \omega)$ will denote the induced labeled sub-poset on $S$.
We say that $S$ 
is {\it $\omega$-natural} if 
the restriction $\omega: S\to \bN^+$ 
is order-preserving. 
We denote by $\mathcal N_\omega(P)$ the set of
all $\omega$-natural subsets of $V(P)$ and 
$\mathcal I_\omega (P)$ the set of all 
$\omega$-natural ideals of $(P, \omega)$. 
For convenience, 
we will always treat the empty set $\emptyset$ as a 
$\omega$-natural ideal of
$(P, \omega)$. 

Note that, for any finite poset $P$,  
order-preserving labelings and order-reversing labelings for $P$
always exist.

Now, for any labeled poset $(P, \omega)$, we define a directed graph 
$\mathcal G (\Pw)$ as follows.

The set of vertices of $\mathcal G (\Pw)$ is the set 
of ideals of $P$, i.e. elements of $\mathcal I (P)$. Two elements 
$I, J\in \cI (P)$ are connected by an arrow pointing to $I$ if and only if
$I\subsetneq J$ and $J\backslash I\in \mathcal N_\omega (P)$. 

We call the directed graph $\mathcal G (\Pw)$
the {\it $\omega$-graph} of the labeled poset $(P, \omega)$. 
We let $A (P, \omega)$ denote the adjacency matrix of the $\omega$-graph 
$\mathcal G (\Pw)$.
Note that we can always arrange the ideals of $P$ so that 
the matrix $A (P, \omega)$ is strictly
upper-triangular. In particular, the matrix  $A (P, \omega)$
is always nilpotent. In this paper 
we will always assume that 
the matrix $A (P, \omega)$ 
is strictly upper-triangular.

Recall that a (directed) {\it path} of length $k$ in a directed graph $D$ is a list
$v_0,\ldots,v_k$ of vertices of $D$ such that $(v_i,v_{i+1})$ is an arc 
(i.e., an edge directed toward to $v_{i+1}$)
in $D$ for $0 \leq i<k$.  Define a (directed) {\it multi-path} of length $k$ in $D$
to be a list $v_0,\ldots,v_k$ of vertices such that, for $0 \leq i<k$,
either $(v_i,v_{i+1})$ is an arc
in $D$ or $v_{i+1}=v_i$.
For each $k\in \bN^+$,
we denote by ${\mathcal A}_k (P, \omega)$ the set of all directed
multi-paths of length $k$ from
$\mathcal G (\Pw)$ connecting the elements $\emptyset$ and $P$
and ${\mathcal C}_k (P, \omega)$ the set of all directed
paths of length $k$ of
$\mathcal G (\Pw)$ connecting the elements $\emptyset$ and $P$. We let
$a_k(P, \omega)$ and $c_k(P, \omega)$
the cardinalities of 
${\mathcal A}_k(P, \omega)$ and ${\mathcal C}_k(P, \omega)$, respectively.

\subsection{Order Polynomials of Labeled Posets}\label{S2.1}

Let $(P, \omega)$ be a labeled poset.
Recall that a map $f: P\to [n]$ with $n\in \bN^+$ is said to be
{\it $\omega$-order-preserving}  if $f$ is order-preserving 
and, for any $x>y$ in $P$ with 
$\omega (x) < \omega (y)$, we have $f(x)>f(y)$. 

The following two lemmas follow immediately from the definitions.

\begin{lemma}\label{L2.1}
Let $(P, \omega)$ be a labeled poset.

$(a)$ If $\omega$ 
is a natural labeling, i.e. $\omega: P\to \bN^+$ is order-preserving, then
a map $\varphi: P\to [n]$ is 
$\omega$-order-preserving if any only of it is 
order-preserving. In this case, any subset of $P$ is 
$\omega$-natural.

$(b)$ If $\omega: P\to \bN^+$ is order-reversing, then 
a map $\varphi: P\to [n]$ is 
$\omega$-order-preserving if and only if it is 
strictly order-preserving.
A subset $S\subset P$ is 
$\omega$-natural 
if and only of $S$ with induced poset structure from $P$ is 
an anti-chain.
\end{lemma}

\begin{lemma}\label{L2.2}
An order-preserving map $f: P\to [n]$ is  
$\omega$-order-preserving if and only if 
$f^{-1}(k)$ is  $\omega$-natural  for any $k\in [n]$.
\end{lemma}

\begin{lemma}\label{L2.3}
For any $n\in \bN^+$, 
the set of $\omega$-order-preserving maps  $f: P\to [n]$ is in one-to-one
correspondence with the set $\mathcal A_n(P, \omega)$. Hence, the number
of $\omega$-order-preserving maps  $f: P\to [n]$ equals $a_n(\Pw)$.
\end{lemma}
\pf
For any $\omega$-order-preserving map  $f: P\to [n]$, 
Set $I_k=\{x\in P| f(x)\leq k\}$ for any $1\leq k\leq n$. 
It is easy to see that $\emptyset\prec I_1 \prec \cdots \prec I_n=P$
is a multi-path of the directed graph  $\mathcal G(P, \omega)$. 
Conversely, for any  multi-path  
$\emptyset\prec I_1 \prec \cdots \prec I_n=P$  of $\mathcal G(P, \omega)$,
we define a map  $f: P\to [n]$ by setting 
$f(I_1)=1$ and $f(I_{k}\backslash I_{k-1})=k$ 
for any $2\leq k\leq n$. It is easy to check that
the map $f$ is $\omega$-order-preserving and
the correspondences defined above are inverse to each other.
\epfv

Now, for any labeled poset $(P, \omega)$,
we define a function $\theta (P, \omega ; \, \cdot\, ): \bN^+ \to \bN$ by setting 
$\theta(P, \omega; n)$ $(n\in \bN^+)$ to be the number of 
$\omega$-order-preserving maps $f: P \to [n]$.  
From Lemma \ref{L2.3}, we see that 
$\theta (\Pw; n)=a_n(P, \omega)$ 
for any $n \in \bN^+$. 

Next we give a proof of the following well-known theorem.

\begin{theo-def}\label{Theo-Def}
For any labeled poset $(P, \omega)$, there exists a unique polynomial
$\Omega (\Pw; t)$ such that $\Omega (\Pw, n)=\theta (\Pw; n)$
for any $n\in \bN^+$.
\end{theo-def}

Note that $\Omega (\emptyset, t)=1$.

To prove this theorem, we first define 
the following $|\mathcal I(P)| \times |\mathcal I(P)|$ matrices.

\begin{align}
\Phi (P, \omega) & :=\ln (1+A (P, \omega))=\sum_{n=1}^\infty 
\frac {(-1)^n A^{n-1}(P, \omega)} n, \label{DefForPhi} \\
\Theta (P, \omega; t)& :=e^{t\Phi(\Pw)}=\sum_{n=0}^\infty 
\frac {\Phi^n(P, \omega)}{n!} t^n, \label{DefForTheta} 
\end{align} 
where $t$ is a formal variable.

First note that $\Phi (P, \omega)$
is a strict upper triangular matrix, for 
the adjacency matrix 
$A (P, \omega)$ is.
In particular, 
$\Phi (P, \omega)$ is nilpotent. 
Hence $\Theta (\Pw, t)$ is a upper triangular matrix
with $1$ on the diagonal and each entry 
in the polynomial algebra $\bQ[t]$.

The theorem above follows immediately from the following lemma.

\begin{lemma}\label{L2.5}
For any  $I, I'\in \cI (P)$, let $\theta (I, I'; t)$
denote  the $(I, I')^{th}$ entry of the matrix
$\Theta (P, \omega; t)$. Then, 
\begin{align}\label{EqFor-theta}
\theta (I, I'; t)=\begin{cases} 0  &\text{\quad   if $I\nsubseteq I'$;}\\ 
\Omega (I'\backslash I, \omega, t) & \text{\quad  if $I\subseteq I'$.}
\end{cases}
\end{align}
In particular, 
$\theta (\emptyset, P; t)$ is the same as
the order polynomial 
$\Omega (\Pw; t)$.
\end{lemma}

\pf
For any $n\in \bN^+$, we have
\BQn
\Theta (P, \omega; n)&=& e^{n\Phi (P, \omega)}\\
&=&e^{n\ln (1+A(P, \omega) )}\\
&=&(1+A(P, \omega))^n.
\EQn

By a standard result in algebraic graph theory (see \cite{B}), 
the
$(I, I')^{th}$ entry of the matrix 
$(1+A(P, \omega))^n$ equals to the number of directed
multi-paths in the $\omega$-graph $\cG(P, \omega)$ 
connecting $I$ and $I'$. It is easy to see that it is also 
same as 
the number of the directed multi-paths in the $\omega$-graph 
$\cG (I'\backslash I, \omega)$ 
of the induced labeled sub-poset $(I'\backslash I, \omega)$ 
connecting $\emptyset$ and $I'\backslash I$ when $I\subseteq I'$.
Therefore, for any $n\in \bN^+$, we have
$\theta (I, I'; n)=0$ if $I \nsubseteq I'$ and 
$\theta (I, I'; n)=\theta (I' \backslash I, \omega; n)$ if $I \nsubseteq I'$. 
By Lemma \ref{L2.3} 
and the definition of order polynomials, 
we see that Eq.\,(\ref{EqFor-theta}) holds.
\epfv

Next, we derive some properties of order polynomials of labeled
posets $(P, \omega)$ and
also study the combinatorial interpretation of 
the numerical invariants defined by the matrix $\Phi(P, \omega)$.

\begin{theo}\label{s+t}
For any labeled poset $(P, \omega)$ and $s, t\in \bC$,  
we have
\BQ \label{Eq(s+t)}
\Omega (P, \omega; s+t) = \sum_{S\in \cI (P)}
\Omega (S, \omega; s) \Omega (P\backslash S, \omega; t).
\EQ
\end{theo}

\pf First, from Eq.(\ref{DefForPhi}) we have
\BQ
\Theta (P, \omega; s+t)=\Theta (P, \omega; s)\Theta (P, \omega; t).
\EQ

By Lemma \ref{L2.3} and the equation above, we have
\BQn
\Omega (P, \omega; s+t) &=& \theta (\emptyset, P, s+t)\\
&=& \sum_{S\in \cI (P)}
\theta (\emptyset, S; s)\theta (S, P; t)\\
 &=& \sum_{S\in \cI (P)}
\Omega (S, \omega; s)\Omega (P\backslash S, \omega; t).
\EQn
\epfv

The next theorem is a direct consequence of Theorem \ref{s+t} above, 
but is more convenient and useful. 
It not only provides an algorithm for the calculation of 
order polynomials of labeled posets, but also is the main motivation
for the further generalizations that will be discussed in 
Section \ref{S3}. 

First we denote by $\Delta$ the difference operator on the polynomial
algebra $\bC[t]$, i.e. 
\BQn
(\Delta f)(t)=f(t+1)-f(t).
\EQn
for any $f\in \bC[t]$. We also define the linear operator
$\Delta^{-1}: \bC[t]\to \bC[t]$ by setting $(\Delta^{-1} f)(t)$ 
to be the unique polynomial $g(t)$ such that 
\BQn 
g(0)&=&0,\\
(\Delta g)(t) &=& f(t).
\EQn 

\begin{theo}\label{Delta}
For any labeled poset $(P, \omega)$ with $|P|=p\geq 2$, we have
\BQ\label{DeltaEqR}
 \Omega (\Pw; t)=\Delta^{-1}
\sum_{\emptyset \neq S\in \mathcal I_\omega (P)}
\Omega (P\backslash S, \omega; t).
\EQ
\end{theo}

Note that the sum in Eq.\,(\ref{DeltaEqR}) 
runs over a (usually) smaller set, namely the set
$\cI_\omega (P)$ of $\omega$-natural ideals of $(\Pw)$, 
than the one in Eq.\,(\ref{Eq(s+t)}), which runs over the set 
$\cI(P)$ of all ideals of $P$.

\pf From Lemma \ref{L2.2}, for any labeled poset $(Q, \eta)$, 
we have
\BQn
\Omega (Q, \eta; 1)=\begin{cases} 1 \quad &\text{if $(Q, \eta)$ is natural;}\\
 0 \quad &\text{otherwise}.\end{cases}
\EQn
This and 
Eq.\,(\ref{Eq(s+t)}),
setting $s=1$, yield
\BQ\label{DeltaEqR-prime}
\Delta \Omega (\Pw; t)=
\sum_{\emptyset \neq S\in \mathcal I_\omega (P)}
\Omega (P\backslash S, \omega; t).
\EQ
Since $\Omega (\Pw; t)$ has no constant term, 
Eq.\,(\ref{DeltaEqR}) follows immediately from the equation above.
\epfv
  
Finally, to end this subsection,
let us look at a combinatorial interpretation of 
coefficients of order polynomials of labeled posets.

First, from Eq.\,(\ref{DefForTheta}), we have 
\BQ\label{ddt}
\frac {d}{d t} \Theta (P, \omega; 0)=\Phi (P, \omega).
\EQ 

For any $I, I'\in \cI(P, \omega)$,
we define $\phi_{I', I}\in \bQ$ by writing 
$\Phi (P, \omega)=(\phi_{I', I})_{I, I'\in \cI (P)}$. 
From Lemma \ref{L2.5} and Eq.\,(\ref{ddt}), we see that
$\phi_{I', I}$ is equal to
the coefficient of $t$ in the order polynomial 
$\Omega (I'\backslash I, \omega; t)$ if $I\subsetneq I'$ 
and $0$ otherwise. In particular, the quantity $\phi_{I', I}$
depends only on the induced
labeled sub-poset structure on the subset 
$I'\backslash I$ from $(\Pw)$, not on 
the ambient labeled poset $(\Pw)$. This also follows from 
Eq.\,(\ref{EqFor-phi}) below. 
So, from now on,  
$\phi_{I', I}$ will also be denoted 
by $\phi (I'\backslash I; \omega)$ 
when $I\subsetneq I'$. Hence, for any labeled poset $(\Pw)$, we have 
a unique well defined rational number $\phi(\Pw)$.

From Eq.\,(\ref{DefForPhi}), it is easy to see that
the following lemma is true.

\begin{lemma}\label{L2.8}
For any  labeled poset $(P, \omega)$ 
and $I,  I'\in \cI (P, \omega)$, we have
\BQ\label{EqFor-phi}
\phi_{I, I'} = 
\begin{cases}  \sum_{k=1}^{|I'\backslash I|} (-1)^{k-1}
\frac {c_k (I'\backslash I, \omega)}k \quad & \text{if $I\subsetneq I'$;}\\ 
 0 \quad & \text{otherwise.}\end{cases} 
\EQ
In particular, for any non-empty labeled poset $(\Pw)$, we have
\BQ\label{EqFor-phi-2}
\phi(\Pw)=\sum_{k=1}^{|P|} (-1)^{k-1}
\frac {c_k (P, \omega)}k.
\EQ
\end{lemma}

Note that
$c_k (I'\backslash I, \omega)$ is the number of directed
chains of the $\omega$-graph $\cG(\Pw)$ connecting $I$ and $I'$. It is
also equal to 
the number of directed
chains of the $\omega$-graph $\cG(I'\backslash I, \omega)$ of the 
induced labeled sub-poset 
$(I'\backslash I, \omega)$ of $(\Pw)$
connecting $\emptyset$ and $I'\backslash I$.

Interestingly, the sequence of rational numbers $\{ \phi (\Pw) \}$ 
indexed by  the set of labeled posets satisfies 
the following recursion formula. 

\begin{propo}
$(1)$ $\phi (\emptyset) =0$ and $\phi (\bullet)=1$.

$(2)$
For any non-empty poset $P$, set
$n_{(\Pw)}$ equal to $1$ if $P$ is $\omega$-natural and
$0$ otherwise.
Then

\BQ\label{RecurEqForphi}
&{}& \phi (\Pw) =n_{(\Pw)}\\
&{}& \quad 
-\sum_{r=2}^{|P|}\frac {1}{r!} 
\sum_{\emptyset \neq I_1 \vartriangleleft I_2\vartriangleleft \cdots \vartriangleleft I_r=P}
 \phi (I_1, \omega) 
\phi (I_2\backslash I_1, \omega)
\cdots \phi (I_r\backslash I_{r-1}, \omega),\nno
\EQ
 
or equivalently,

\BQ\label{SumPsi=0}
&{}& \sum_{r=1}^{|P|} \frac {1}{r!} \sum_{\emptyset \neq
I_1 \vartriangleleft I_2\vartriangleleft\cdots \vartriangleleft I_r=P}
\phi (I_1, \omega) 
\phi (I_2\backslash I_1, \omega)
\cdots \phi (I_r\backslash I_{r-1}, \omega)\\
&{}& \quad \quad
=n_{(\Pw)}. \nno
\EQ
\end{propo}
\pf $(1)$ is obvious from the definition of $\phi$ and 
Eq.\,(\ref{EqFor-phi-2}). To prove $(2)$, 
first from Eq.\,(\ref{DefForPhi}), we have
\BQn
e^{\Phi (P, \omega)}=1+A(\Pw).
\EQn
Hence we have,
\BQn
\sum_{r=1}^\infty \frac {1}{r!} \Phi^r(P, \omega)=A(\Pw).
\EQn
Note that $\Phi^r(P, \omega)=0$ for any $r>|P|$, since there are 
no 
paths of length greater than $|P|$ in the $\omega$-graph $\cG(\Pw)$. 
Hence, we have
\BQn
\sum_{r=1}^{|P|} \frac {1}{r!} \sum_{
I_1, I_2, \cdots,  I_r=P\in \cI (P) }
 \phi (\emptyset, I_1) \phi (I_1, I_2)
\cdots \phi (I_{r-1}, P)=n_{(\Pw)}.
\EQn

Then, Eq.\,(\ref{SumPsi=0}) follows from 
$(1)$ and Lemma \ref{L2.8}.
\epfv

From Eq.\,(\ref{DefForTheta}) and Lemma \ref{L2.5}, it is easy to see
that we can write down the order polynomials 
$\Omega (\Pw; t)$ in terms of 
the rational numbers $\phi$ as follows.

\begin{propo}\label{P2.10}
For any non-empty
labeled poset $(\Pw)$, we have
\BQn
 &{}&\quad \quad \Omega (\Pw; t)= \\ 
  &{}& 
\sum_{r=1}^{|P|} \left ( \frac 1{r!} 
\sum_{\emptyset \neq I_1\vartriangleleft  I_2
 \vartriangleleft \cdots \vartriangleleft I_r=P}
\phi (I_1, \omega) 
\phi (I_2\backslash I_1, \omega)
\cdots \phi (I_r\backslash I_{r-1}, \omega)
\right ) t^r. 
\EQn
\end{propo}

In terms of the $\phi$'s, we also have the following recursion formula for 
order polynomials of labeled posets.

\begin{propo}
For any finite poset $P$, we have 
\BQ
\Omega' (\Pw; t)&=&\sum_{S \in \mathcal J (P) }
\phi (S, \omega) \Omega (P\backslash S, \omega; t) \\
&=&
\sum_{S \in \mathcal J (P) }
\phi (P\backslash S, \omega) \Omega (S, \omega;  t),\nno
\EQ
or in other words,
\BQ
\Omega' (\Pw; t)&=&\sum_{S \in \mathcal J (P) }
\Omega' (S, \omega; 0) \Omega (P\backslash S, \omega;  t)\\
&=&\sum_{S \in \mathcal J (P) }
\Omega'(P\backslash S, \omega;  0) \Omega' (S, \omega; t).\nno
\EQ
\end{propo}

\pf By applying $\frac {d}{d t}$ to Eq.\,(\ref{DefForTheta}), 
we get
\BQ
\frac {d}{d t} \Theta (\Pw; t)=\Phi (\Pw) \Theta (\Pw; t)
=\Theta (\Pw; t)\Phi (\Pw). 
\EQ
Then the equations in the theorem follow immediately 
from the equation above and Lemma \ref{L2.5}.
\epfv

\subsection{Eulerian Polynomials of Labeled Posets} \label{S2.2}

In this subsection, we give a formula for the Eulerian polynomial 
$e (P, \omega; \lambda)$
of a labeled poset $(\Pw)$ 
in terms of its adjacency matrix $A(\Pw)$. We also derive certain 
 recursion formulas for the Eulerian polynomial 
$e (P, \omega; \lambda)$.

Let $(\Pw)$ be a labeled poset.
Recall that the Eulerian polynomial 
$e (P, \omega; \lambda)$ of $(\Pw)$
is defined by the equation
\BQ\label{Euler-Poly}
 \sum_{n=0}^\infty \Omega (P, \omega; n) \lambda^n=
\frac {e (P, \omega; \lambda )}{(1-\lambda)^{|P|+1}}.
\EQ

\begin{rmk}
It is a standard fact in combinatorics that $e (P, \omega; \lambda)$ 
defined above is indeed a polynomial; however, this will be clear 
from Eq. $(\ref{Formu-e})$ below.
\end{rmk}

We set
\BQ\label{e-tilde}
\tilde e (P, \omega; \lambda )=\frac {e (P, \omega; \lambda )}{(1-\lambda)^{|P|+1}}
\EQ
and
\BQ\label{DefForE}
E(P,\omega; \lambda):=\sum_{n=0}^\infty \Theta (P, \omega; n) \lambda^n.
\EQ
and, for any
$I, I'\in \cI (P)$,
define 
$\tilde e (I, I'; \lambda )$ by writing 
$E(P,\omega; \lambda)=\left (
\tilde e (I, I'; \lambda )\right )_{I, I'\in \cI (P)}$. 

From Lemma \ref{L2.5} and Eq.\,(\ref{Euler-Poly}) and (\ref{DefForE}), 
it is easy to see that we have 
the following lemma.

\begin{lemma}\label{L2.12}
For any
$I, I'\in \cI (P)$, we have
\BQ\label{EqForTilde-e}
\tilde e (I, I'; \lambda )=\begin{cases}
 \tilde e (I'\backslash I,\omega; \lambda )
 \quad &\text{if $I\subseteq I'$;} \\
0\quad &\text{otherwise.} \end{cases}
\EQ
\end{lemma}

In particular, $\tilde e (I, I'; \lambda )$
depends only on the reduced
labeled sub-poset structure on the subset 
$I'\backslash I$ of $P$, not on 
the ambient labeled poset $P$. 
So, from now on,  we will also denote $\tilde e (I, I'; \lambda )$
by $\tilde e (I'\backslash I, \omega; \lambda )$ 
when $I\subseteq I'$.

\begin{theo}\label{TheoFor-E}
For any labeled poset $(\Pw)$, we set $\wtilde A(\Pw):=\text{I} +A(\Pw)$.
Then, we have
\BQ\label{EqFor-E-1}
E(P, \omega; \lambda)=(I-\lambda  \wtilde A(\Pw) )^{-1}
\EQ
and
\BQ\label{EqFor-E-2}
E(P, \omega; \lambda) =
\frac 1{1-\lambda}
\left (I-\frac \lambda{1-\lambda} A(\Pw) \right )^{-1}.
\EQ
\end{theo}

\pf It is easy to see that 
Eq.\,(\ref{EqFor-E-2}) is a direct consequence of Eq.\,(\ref{EqFor-E-1}).
To prove Eq.\,(\ref{EqFor-E-1}),
first
by Eq.\,(\ref{DefForTheta}) and (\ref{DefForE}), we have 
\BQn
E(P, \omega; \lambda)&=& \sum_{n=0}^\infty \Theta (\Pw; n)
\lambda^n \\
&=& \sum_{n=0}^\infty (\text{I}+A(\Pw) )^n \lambda^n\\
&=& \sum_{n=0}^\infty  \wtilde A(\Pw) ^n \lambda^n \\
&=&(I-\lambda  \wtilde A(\Pw) )^{-1}.
\EQn
\epfv

\begin{corol}\label{C2.14}
For any non-empty labeled poset $(\Pw)$, we have
\BQ\label{Formu-tilde-e}
\tilde e(\Pw; \lambda) &=&
\sum_{k=1}^{|P|} a_k(\Pw) \lambda ^k \\
&=&\frac 1{1-\lambda}
\sum_{k=1}^{|P|} c_k(\Pw) \left (\frac \lambda {1-\lambda}\right )^{k} \nno
\EQ
and
\BQ\label{Formu-e}
e(\Pw; \lambda)=
\sum_{k=1}^{|P|} c_k(\Pw) \lambda^k (1-\lambda)^{|P|-k}. 
\EQ
\end{corol}
\pf Eq.\,(\ref{Formu-tilde-e}) follows from  Eq.\,(\ref{EqFor-E-1}), 
(\ref{EqFor-E-2}) and (\ref{EqForTilde-e});
 while  Eq.\,(\ref{Formu-e}) follows from 
Eq.\,(\ref{Formu-tilde-e}) and (\ref{e-tilde}). 
\epfv.

Recall that, for any non-empty directed graph $G$, the 
{\it chain polynomial} $c(G; \mu)$ of $G$ is defined to be
\BQn
c(G, \mu)=\sum_{k=0}^{\infty} c_k(G)\mu^k=\sum_{k=0}^{|G|} c_k(G) \mu^k,
\EQn
where $c_0(G)=1$ and $c_k(G)$ $(k\geq 1)$ is the number of 
directed chains of length $k$ in $G$. 

For any non-empty labeled poset $(P, \omega)$, we set
$G(P, \omega)$ to be the directed graph 
${\mathcal G}(P, \omega)\backslash \{\emptyset, P\}$.
Then,  
for any $k\geq 1$, it is easy to see that $c_k(P, \omega)$ 
is the same as 
the number $c_{k-1}(G(P, \omega))$. Hence, 
by Eq.\,(\ref{Formu-tilde-e}) and (\ref{e-tilde}), we have 
the following identities.

\begin{propo}
For any non-empty labeled poset $(\Pw)$, we have
\BQn
\tilde e(\Pw; \lambda)&=& \frac \lambda {(1-\lambda)^2} 
c(G; \frac \lambda{1-\lambda})\\
e(\Pw; \lambda)&=&\lambda (1-\lambda)^{|P|-1} c(G; \frac \lambda{1-\lambda}).
\EQn
\end{propo} 

\begin{theo}\label{T2.15}
Let $(\Pw)$ be a non-empty
labeled poset. Then we have 
\BQ\label{recurFor-e}
e(\Pw; \lambda) =\lambda \sum_{\emptyset \neq S\in \cI_\omega (P)} 
(1-\lambda)^{|S|-1}e (P\backslash S, \omega; \lambda)
\EQ
and
\BQ\label{recurFor-tilde-e}
\tilde e(\Pw; \lambda)=\frac \lambda {1-\lambda} 
 \sum_{\emptyset \neq S\in \cI_\omega (P)} 
\tilde e(P\backslash S, \omega; \lambda).
\EQ
\end{theo}
It is easy to see that  Eq.\,(\ref{recurFor-e}) and 
(\ref{recurFor-tilde-e}) are equivalent to each other, 
so we need to prove only one of them.

{\it First Proof:} 
First, from the equation
\BQn
\lambda (I-\lambda \wtilde A(\Pw))^{-1}\wtilde A(\Pw)
=(I-\lambda \wtilde A(\Pw))^{-1}-I
\EQn
and Eq.\,(\ref{EqFor-E-1}), we get
\BQn
\lambda \wtilde A(\Pw) E(\Pw; \lambda)= E(\Pw; \lambda)-I.
\EQn

Now by comparing the $(\emptyset, P)^{th}$ entries of 
both sides of 
the equation above, we have
\BQn
\lambda \sum_{S\in \cI_\omega (P)} 
\tilde e (P\backslash S, \omega; \lambda)&=&
\lambda \tilde e(\Pw, \lambda)+
\lambda \sum_{\emptyset \neq S\in \cI_\omega (P)} 
\tilde e (P\backslash S, \omega; \lambda)\\
&=&\tilde e(\Pw, \lambda).
\EQn
By solving for $\tilde e(\Pw, \lambda)$ in the equation above, 
we get Eq.\,(\ref{recurFor-tilde-e}).
\epfv

We also can prove  Eq.\,(\ref{recurFor-e}) directly as follows.

{\it Second Proof:} 
First, by Theorem \ref{Delta}, we have
\BQ
\Omega (\Pw; t+1)-\Omega (\Pw; t)
=\sum_{\emptyset \neq S\in \cI_\omega (P)} 
\Omega (P\backslash S, \omega; t).
\EQ
Hence we have
\BQn
\sum_{n=0}^\infty (\Omega (\Pw; n+1)-\Omega (\Pw; n) )\lambda^n 
=\sum_{\emptyset \neq S\in \cI_\omega (P)} \sum_{n=0}^\infty  
\Omega (P\backslash S, \omega; n)\lambda^n. 
\EQn
and since $\Omega (\Pw; 0)=0$ (see Proposition \ref{P2.10}), we have
\BQn
\frac {1-\lambda}{\lambda} \sum_{n=0}^\infty \Omega (\Pw; n)  \lambda^n 
=\sum_{\emptyset \neq S\in \cI_\omega (P)} 
\frac {e (P\backslash S, \omega; \lambda)}{(1-\lambda)^{|P\backslash S|+1}},
\EQn
and
\BQn
\frac {1-\lambda}{\lambda}\frac {e (\Pw; \lambda)}{(1-\lambda)^{|P|+1}}
=\sum_{\emptyset \neq S\subset \cI_\omega (P)} 
\frac {e (P\backslash S, \omega; \lambda)}{(1-\lambda)^{|P\backslash S|+1}}.
\EQn
Hence we have Eq.\,(\ref{recurFor-tilde-e}).
\epfv

Finally, let us apply our recursion formula (\ref{recurFor-tilde-e}) 
to the anti-chains $A_n$ 
$(n\geq 1)$. Note that, for any fixed $n\in \bN^+$,
the labeled poset $(A_n, \omega)$ is always 
a naturally labeled poset for any labeling $\omega$ and hence, for fixed $n$, all
$(A_n,\omega)$ are equivalent.
The Eulerian polynomial 
$A_n(\lambda)=e(A_n, \omega; \lambda)$ 
of $(A_n, \omega)$ is known as the $n^{th}$
{\it Eulerian polynomial}.

\begin{corol}
For any $n\geq 2$, we have
\BQ
A_n(\lambda)=\lambda \sum_{k=1}^n \binom nk A_{n-k}(\lambda)(1-\lambda)^{k-1}. 
\EQ
\end{corol}
\pf This is a direct consequence of Eq.\,(\ref{recurFor-tilde-e}), 
because any subset
$S$ of $(A_n; \omega)$ is in $\cI_\omega (A_n)$ and 
$A_n\backslash S=A_{n-|S|}$.
\epfv
 
Actually, we have another recursion formula for $A_n(\lambda)$ $(n\in \bN^+)$, 
which is better in the sense that it involves only $A_{n-1}(\lambda)$.

\begin{propo}
For any $n\geq 1$, we have
\BQ\label{A-Prime}
A_{n}(\lambda)=\lambda (1-\lambda)A'_{n-1}(\lambda)+n\lambda A_{n-1}(\lambda).
\EQ
\end{propo}

\pf Note that, for any $n\geq 1$,  we have
$\Omega (A_n, \omega;  t)=t^n$. 
(This follows from the definition of the order polynomial.)  
Therefore
\BQ
\frac {A_n(\lambda)}{(1-\lambda)^{n+1}}=\sum_{k=0}^\infty k^n \lambda^k 
=\left (\lambda \frac d{d \lambda}\right )^n \frac 1{1-\lambda}.
\EQ
Hence we have 
\BQn
\frac {A_n(\lambda)}{(1-\lambda)^{n+1}}=\lambda \frac d{d \lambda}
\frac {A_{n-1} (\lambda)}{(1-\lambda)^{n}}=
\frac{\lambda (1-\lambda)A'_{n-1}(\lambda)+n\lambda A_{n-1}(\lambda)}{(1-\lambda)^{n+1}}
\EQn
and Eq.\,(\ref{A-Prime}) follows immediately.
\epfv

\subsection{Applications to  
Bernoulli Numbers and Bernoulli Polynomials}\label{S2.3}

In this subsection, we use the results on order polynomials of labeled posets 
derived in Subsection \ref{S2.2} to give some combinatorial results on Bernoulli numbers 
and  Bernoulli polynomials.

Consider strictly labeled 
shrubs $(S_n, \omega)$ ($n\in \bN$), i.e.  the labeling
$\omega: S_n\to \bN^+$ is an order-reversing map. 
Note that, by Lemma \ref{L2.1}, the order polynomial 
$\Omega (S_n, \omega; t)$ is 
same as the strict order polynomial 
$\bar \Omega (S_n; t)$ of the unlabeled poset $S_n$. Also note that, 
the only $\omega$-natural ideals of $(S_n, \omega)$ are
$\emptyset$ and the set consisting of the unique minimum element of $S_n$. 

First let us recall that the Bernoulli numbers $b_n$ ($n\in \bN$) 
 are defined by the generating function
\BQ\label{GeneFctn-bn}
\frac x{e^x-1}=\sum_{n=0}^\infty b_n \frac {x^n}{n!}
\EQ
and the Bernoulli polynomials
$B_n(t)$ ($n\in \bN$) 
are  
defined by 
\BQ \label{GeneFctn-Bn}
\frac {xe^{tx}}{e^x-1}=\sum_{n=0}^\infty B_n(t) \frac {x^n}{n!}.
\EQ 

\begin{propo}\cite{WZ}
For the labeled poset $(S_n, \omega)$ $(n\in \bN^+)$, we have
\BQ \label{E2.3.1}
\phi (S_n, \omega)=b_n
\EQ
and
\BQ\label{E2.3.2}
 \Omega (S_n, \omega; t)=\int_0^t B_n(u)\, du.
\EQ
\end{propo}

Two different proofs are given in \cite{WZ}. 
Here we give yet another one.

\pf Clearly, we need to prove only Eq.\,(\ref{E2.3.2}), for $\phi (S_n, \omega)$ and 
$b_n$ are the coefficients of $t$ in polynomials 
$\Omega (S_n, \omega; t)$ and $\int_0^t B_n(u)\, du$, respectively.
By applying Eq.\,(\ref{DeltaEqR}), we have
\BQn
\Delta \Omega (S_n, \omega; t)= \Omega (A_n, \omega; t)=t^n.
\EQn

Since both $\Omega (S_n, \omega; t)$ and $\int_0^t B_n(u)\, du$
have no constant terms, it suffices to show that 
$\int_0^t B_n(u)\, du$ also satisfies the equation above.

First, from the well-known equation $\Delta B_n(u)=nu^{n-1}$, 
we have
\BQn
 \int_0^t \Delta  B_n(u)\, du=t^n.
\EQn
But
\BQn
 \int_0^t \Delta  B_n(u)\, du &=& \int_0^t B_n(u+1)\, du
-\int_0^t  B_n(u)\, du\\
&=&
\int_1^{t+1}   B_n(u)\, du-\int_0^t  B_n(u)\, du\\
&=& \Delta\int_0^t  B_n(u)\, du-\int_0^1  B_n(u)\, du.
\EQn
Since $\int_0^1 \frac{xe^{tx}}{e^x-1}\, dt =1$,  we have
$\int_0^1  B_n(u)\, du=0$ for any $n\in \bN^+$. Therefore, we see that
\BQn
\Delta\int_0^t  B_n(u)\, du=t^n.
\EQn
Hence, we have Eq.\,(\ref{E2.3.2}). \epfv

Note that Eq.\,(\ref{E2.3.2}) provides a combinatorial interpretation for
the Bernoulli polynomials $B_n(t)$ $(n\in \bN^+)$. To consider the 
combinatorial interpretation of 
the Bernoulli numbers $b_n$ $(n\in \bN^+)$, first, 
from Eq.\,(\ref{E2.3.1}) and 
(\ref{EqFor-phi-2}), 
we have
\BQ
b_n=\sum_{k=1}^{n+1} (-1)^{k-1}
\frac {c_k (S_n, \omega)}k.
\EQ
Therefore, $b_n$ is just the weighted alternating sum 
of the numbers of directed paths of the $\omega$-graph ${\mathcal G}(S_n, \omega)$
of the labeled poset
$(S_n, \omega)$ connecting $\emptyset$ and $S_n$, 
where such a directed path of length $k$ $(k\in \bN^+)$
is weighted by $\frac 1k$.
By counting the numbers $c_k(S_n, \omega)$ $(k\in \bN^+)$ 
of directed chains in $\cG (S_n, \omega)$ connecting $\emptyset$ and $S_n$,
it is easy to see that
\BQ
c_1 (S_n, \omega)&=& 0; \\
c_k (S_n, \omega)&=&\sum_{\substack{n_1+\cdots +n_{k-1}=n\\ n_i\geq 1}} 
\binom{n}{n_1,  \cdots, n_{k-1}}
\EQ
for any $k\geq 2$. 
Therefore, we have

\BQ \label{PhiOfSn}
b_n &=&\sum_{k=1}^n (-1)^{k} \frac {1}{k+1}
\sum_{\substack{n_1+\cdots +n_{k}=n\\ n_i\geq 1}} 
\binom{n}{n_1,  \cdots, n_{k}}.
\EQ

On the other hand, 
from the generating function Eq.\,(\ref{GeneFctn-Bn}),
it is easy to derive that
\BQ \label{bn}
\quad \quad \quad b_n &=&\sum_{r=1}^n (-1)^{r} n! \sum_{\substack{n_1+\cdots +n_r=n\\ n_i\geq 1}}
\frac 1{(n_1+1)!}\frac 1{(n_2+1)!}\cdots \frac 1{(n_r+1)!}.
\EQ

Hence, we get the following identity.

\begin{lemma}
For any $n\geq 1$, we have
\BQ
&{}& \sum_{k=1}^n (-1)^{k}  \frac {1}{k+1}
\sum_{\substack{n_1+\cdots +n_{k}=n\\ n_i\geq 1}} 
\frac 1{n_1!}\frac 1{n_2!}\cdots \frac 1{n_k!} \\
&{}& \quad =\sum_{k=1}^n (-1)^{k} \sum_{\substack{n_1+\cdots +n_k=n\\ n_i\geq 1}}
\frac 1{(n_1+1)!}\frac 1{(n_2+1)!}\cdots \frac 1{(n_k+1)!}. \nno
\EQ
\end{lemma}

The identity above can also be proved directly by the method 
of generating functions as follows.

{\it 2nd Proof} : First, consider the generating function
\BQ\label{1+U}
\frac {\ln (1+U)}U=\sum_{k=0}^\infty (-1)^{k} \frac {U^k}{k+1}.
\EQ
Now set $U=e^x-1$, the right hand side of the equation above gives
\BQn
&{}& \sum_{k=0}^\infty (-1)^{k} \frac 1{k+1} \left (\sum_{i=1}^\infty \frac {x^i}i \right )^k\\
&=&\sum_{n=0}^\infty \left (
\sum_{k=1}^n (-1)^{k}  \frac {1}{k+1}
\sum_{\substack{n_1+\cdots +n_{k}=n\\ n_i\geq 1}} 
\frac 1{n_1!}\frac 1{n_2!}\cdots \frac 1{n_k!} 
\right ) x^n. 
\EQn

On the other hand, the left hand side of Eq.\,(\ref{1+U}) gives
\BQ
\frac {\ln (1+U)}U= \frac {\ln (1+(e^x-1))}{e^x-1}= \frac x{e^x-1}.
\EQ
Note that 
\BQn
  \frac x{e^x-1} =
 \sum_{n=0}^\infty  \left (
\sum_{k=1}^n (-1)^{k} \sum_{\substack{n_1+\cdots +n_k=n\\ n_i\geq 1}}
\frac 1{(n_1+1)!}\cdots \frac 1{(n_k+1)!}
 \right ) x^n.
\EQn
By comparing the coefficients of $x^n$ in Eq.\,(\ref{1+U}) 
and using the expression
of $b_n$ in Eq.\,(\ref{bn}), 
we are done.
\epfv

One interesting observation is that, from Eq.\,(\ref{PhiOfSn}), it is easy to see that
$(n+1)!b_n =(n+1)!\phi (S_n, \omega)$
 is always a positive integer for any $n\in \bN$. Hence we have the following interesting 
result about Bernoulli numbers.

\begin{corol}
The Bernoulli's number $b_n$ $(n\in \bN)$
can be written as the following {\bf fraction}
\BQn
b_n= \frac {(n+1)! \sum_{k=1}^n (-1)^{k}  \frac {1}{k+1}
\sum_{\substack{n_1+\cdots +n_{k}=n\\ n_i\geq 1}} \binom{n}{n_1, n_2, \cdots, n_{k}}}
{(n+1)!}.
\EQn
\end{corol}

In particular, the denominator of $b_n$ in reduced form 
is always a factor of $(n+1)!$. This fact is not obvious from 
the expression of $b_n$ given by Eq.\,(\ref{bn}).

\renewcommand{\theequation}{\thesection.\arabic{equation}}
\renewcommand{\therema}{\thesection.\arabic{rema}}
\setcounter{equation}{0}
\setcounter{rema}{0}

\section{\bf A New Approach to Order Polynomials of Labeled Posets and 
Their Generalizations}\label{S3}

In this section, motivated by the property of 
order polynomials of labeled posets 
given by Theorem \ref{Delta} and the family of invariants for rooted trees
constructed in \cite{Z}, 
we introduce a family of invariants 
for labeled posets.
We show that order polynomials $\Omega (\Pw; t)$ and  
the invariant $\tilde e(\Pw; \lambda)$ 
(see Eq.\,(\ref{e-tilde})) 
belong to this family of invariants.
We also consider the
known quasi-symmetric function invariant $K(\Pw; x)$ of labeled posets
$(\Pw)$, which suggests
further generalizations of our construction.

\subsection{A Family of Invariants of Labeled Posets}\label{S3.1}

Let $A$ be any vector space over a field $k$ and $\Xi$ a $k$-linear operator
on $A$. We fix one element 
$a\in A$ and define an invariant $\Psi (\Pw)$ for
labeled posets $(\Pw)$ by the following algorithm.

\begin{Algo}\label{Algo}
$(1)$ We first set 
$\Psi (\emptyset)=a$.

$(2)$  For any non-empty labeled poset $(\Pw)$, 
we define  $\Psi (\Pw)$ recursively by
\BQ\label{Algo-eq}
\Psi (\Pw)= \Xi \sum_{\emptyset \neq S\in \cI_\omega (P)} \,  
\Psi (P\backslash S, \omega)
\EQ
\end{Algo}

The following lemma is obvious.

\begin{lemma} \label{L3.1}
 Let  $\Gamma(\Pw)$ be an $A$-valued invariant  
for labeled posets $(\Pw)$.
Then, 
$\Gamma(\Pw)$ can be re-defined and calculated by Algorithm
\ref{Algo} for some $k$-linear map $\Xi$ if and only if 
\begin{enumerate}
\item $\Gamma (\emptyset)=a$, and 
\item $\Gamma$
 satisfies Eq.\,$(\ref{Algo-eq})$ for every 
labeled poset $(\Pw)$.
\end{enumerate}
\end{lemma}

\begin{rmk}
In general,
the restriction of Algorithm \ref{Algo}  
to naturally labeled rooted forests is 
not same as Algorithm  $3.1$ in \cite{Z}, even when 
the vector space has an algebra structure.
 This is because the invariants $\Psi(\Pw)$ defined
by Algorithm \ref{Algo} in general do not satisfy the following  
equation.
\BQn
\Psi(P_1 \sqcup P_2\sqcup \cdots \sqcup P_k)=\prod_{i=1}^k,
\Psi(P_i)
\EQn
where $P_1 \sqcup P_2\sqcup \cdots \sqcup P_k$ is the disjoint union
of naturally labeled posets $P_i$ $(1\leq i\leq k)$. But, 
when $A=\bC[t]$, $a=1$ and $\Xi=\Delta^{-1}$, 
both algorithms give the order polynomials of rooted forests. 
$(\text{See Proposition \ref{Refor-Omega} below.})$
\end{rmk}

\subsection{Re-formulation for order Polynomials 
and Eulerian Polynomials of Labeled Posets}\label{S3.2}
 
From Theorem \ref{Delta} and Lemma \ref{L3.1}, 
it is easy to see that the order polynomials $\Omega (\Pw; t)$ 
of labeled posets $(\Pw)$
can be reformulated as follows.

\begin{propo}\label{Refor-Omega}
 Let $\Psi$ be the 
 invariant 
defined by Algorithm \ref{Algo} with 
 $A=\bQ [t]$, $a=1$ and $\Xi=\Delta^{-1}$.
 Then, for any finite poset $P$, 
we have $\Psi (\Pw)=
\Omega (\Pw; t)$.
\end{propo}

From Theorem \ref{T2.15} and Lemma \ref{L3.1}, 
it is easy to see that the invariant $\tilde e(\Pw; \lambda)$
of labeled posets $(\Pw)$
can be reformulated as follows.

Let $A$ be the localization $\bQ [\lambda]_{(1-\lambda)}$
of the polynomial algebra $\bQ [\lambda]$ at $1-\lambda$, 
i.e. $A=\bQ [\lambda, (1-\lambda)^{-1}]$. Let $M_\lambda$ be
the linear operator of $A$
defined by the multiplication by 
$\frac \lambda {1-\lambda}$. Then we have the following proposition.

\begin{propo}
Let $\Psi$ be the 
invariant defined by Algorithm \ref{Algo} with 
$A=\bQ [\lambda]_{(1-\lambda)}$, 
$a= \frac 1{1-\lambda}$
and $\Xi=M_\lambda$.
Then, for any labeled poset $(\Pw)$, 
we have $\Psi (\Pw)=
\tilde e (\Pw; \lambda)$.
\end{propo}

Since the Eulerian polynomial $e(\Pw; \lambda)=(1-\lambda)^{|P|+1}
\tilde e(\Pw; \lambda)$ for any labeled posets $(\Pw)$, 
we see that Eulerian polynomials $e(\Pw; \lambda)$ 
can  be recovered up to the factor $(1-\lambda)^{|P|+1}$
from an invariant defined Algorithm \ref{Algo}.

\subsection{Quasi-Symmetric Functions $K$ of Labeled Posets}\label{S3.3}

First, let us recall that the following well-known quasi-symmetric 
function invariant $K(\Pw; x)$ for labeled poset $(\Pw)$.

Let $x=(x_1, x_2, \cdots )$ be a sequence of commutative variables and let
$\bC[[x]]$ be 
the formal power series algebra in $\{x_k |k\geq 1\}$ over $\bC$.
For any labeled poset 
$(\Pw)$ and any map $\sigma : P\to \bN^+$ of sets, 
we set
$x^\sigma:=\prod_{i=1}^\infty x_i^{|\sigma^{-1}(i)|}$ and
define
\BQ\label{DefBarK}
 K (\Pw; x)=\sum_{\sigma} x^\sigma, 
\EQ
where the sum runs over the set of all 
$\omega$-order-preserving maps
$\sigma : P\to \bN^+$.

Recall that an element $f\in \bC[[x]]$ is said to be 
{\it quasi-symmetric} if the degree of $f$ is bounded,
and for any $a_1, a_2, \cdots, a_k \in \bN^+$, 
$i_1<i_2<\cdots <i_k$ and
 $j_1<j_2<\cdots <j_k$, 
the coefficient of the monomial 
$x_{i_1}^{a_1}x_{i_2}^{a_2}\cdots x_{i_k}^{a_k}$ is 
always same as the coefficient of the monomial 
$x_{j_1}^{a_1}x_{j_2}^{a_2}\cdots x_{j_k}^{a_k}$. 
 For more general studies on 
quasi-symmetric functions, 
see \cite{G}, \cite{T}, \cite{MR} and \cite{St3}.

From Eq.\,(\ref{DefBarK}), it is easy to check 
that, for any labeled
poset $(\Pw)$,  $K(\Pw; x)$ is quasi-symmetric.

Next, let us derive a 
recursion formula for 
the quasi-symmetric function invariant
$K(\Pw; x)$ of a labeled poset $(\Pw)$.

Let  $S: \bC[[x]]\to \bC[[x]]$ be the shift operator
defined first by setting
\BQn
S (1)&=&1; \\
S (x_m)&=&x_{m+1} \text{\quad for any $m\geq 2$}
\EQn
and then extending it
to the unique $\bC$-algebra
homomorphism from $\bC[[x]]$ to $\bC[[x]]$.
For any $m\in \bN^+$, we define the linear operator 
$\Lambda_m$ of  $\bC[[x]]$ by setting
\BQ\label{Lambda-m}
 \Lambda_m =\sum_{k=1}^\infty x_k^m S^{k} 
\EQ

\begin{propo}\label{Recur-K}
For any non-empty labeled poset $(\Pw)$,  we have
\BQ
 K (\Pw; x) 
=\sum_{\emptyset \neq S\in \cI_\omega (P)} \,  \Lambda_{|S|}
K (P\backslash S, \omega; x).
\EQ
\end{propo}

This proposition is a generalization of Lemma $7.3$ 
in \cite{Z} and the proof is
also similar to the proof there.

\pf
Let $W$ be the set of $\omega$-order-preserving maps 
$\sigma: P\to \bN^+$. For each $k\in \bN^+$, we set 
 $W_k$ to be the set of $\sigma \in W$ such that 
$\text{min} (\sigma (P))=k$.
 Clearly, $W$ is the disjoint union
of the $W_k$ $(k\geq 1)$. Also,
from the definitions of $K(\Pw, x)$ and the shift operator $S$, it is 
easy to see that,   for any labeled poset
$(\Pw)$, we have
\BQ\label{AA}
\sum_{\substack{ \sigma \in W \\ \sigma (P)\subset \bN^+ \backslash [k]}} 
x^\sigma =S^k K(\Pw, x).
\EQ

Note that, for any $\sigma\in W_k$, 
$\sigma^{-1}(k)$ is always 
a non-empty 
$\omega$-natural ideal, 
i.e. $\emptyset \neq \sigma^{-1}(k)\in \cI_\omega (P)$.
Now let $W(P \backslash S, \omega)$ be the set of 
all $\omega$-order-preserving maps $\mu$ from the labeled poset
$(P \backslash S, \omega)$ to $\bN^+$ and 
consider

\begin{align*}
\sum_{\sigma\in W_k} x^\sigma &= \sum_{\emptyset\neq S\in \cI_\omega (P)}
\sum_{\substack{\sigma\in W_k\\
\sigma^{-1}(k)=S}}x^\sigma\\
&= \sum_{\emptyset\neq S\in \cI_\omega (P)}
 x_k^{|S|} \sum_{\substack{\mu \in W(P \backslash S, \omega)\\
\mu (P \backslash S)\subset \bN^+\backslash [k]} } x^\mu \\
&= \sum_{\emptyset\neq S\in \cI_\omega (P)}
  x_k^{|S|} S^k K (P \backslash S, \omega; x), 
\end{align*}
where the last equality follows from Eq.\,(\ref{AA}) and 
the definition of $K(\Pw, x)$ (see Eq.\,(\ref{DefBarK}).).

Therefore, we have
\BQ
K(\Pw; x)&=&\sum_{k=1}^\infty \sum_{\sigma\in W_k} x^\sigma\\
&=&\sum_{k=1}^\infty \sum_{\emptyset\neq S\in \cI_\omega (P)}
x_k^{|S|} S^k K (P \backslash S, \omega; x) \nno\\
&=&\sum_{\emptyset\neq S\in \cI_\omega (P)} 
\sum_{k=1}^\infty x_k^{|S|} S^k K (P \backslash S, \omega; x)\nno\\
&=&\sum_{\emptyset\neq S\in \cI_\omega (P)} \Lambda_{|S|}
K (P \backslash S, \omega; x).\nno
\EQ
\epfv

It does not seem possible to recover $K(\Pw, x)$ from Algorithm \ref{Algo}, 
Rather, the following generalization of  Algorithm \ref{Algo} is suggested.

Let $A$ and $a\in A$ as before. Let $\{ \Xi_m |m\in \bN^+\}$ be
 a sequence of $k$-linear operators of $A$. 
We define an invariant $\Psi (\Pw)$ for
labeled posets $(\Pw)$ 
by the following algorithm.

\begin{Algo}\label{G-Algo}
$(1)$ 
We first set 
$\Psi (\emptyset)=a$.

$(2)$  For any labeled poset $P$ with $|P|\geq 2$, 
we define  $\Psi (\Pw)$ recursively by
\BQ
\Psi (\Pw)=\sum_{\emptyset \neq S \in \cI_\omega (P)} \, 
\Xi_{|S|}  \Psi (P\backslash S, \omega).
\EQ
\end{Algo}

Hence, with $A=\bC [[x]]$, $a=1$ and 
$\Xi_m=\Lambda_m$ $(m\in \bN^+)$, the invariant defined 
by Algorithm \ref{G-Algo} is same as $K(\Pw; \lambda)$.
From Eq.\,(\ref{recurFor-e}) and the fact $e(\emptyset; \lambda)=1$, we see
that the Eulerian polynomials $e (\Pw; \lambda)$ of labeled posets
$(\Pw)$ can also  be calculated by Algorithm \ref{G-Algo} with $A=\bC [\lambda]$, 
$a=1$ and $\Xi_m$ $(m\in \bN^+)$ the linear operators of multiplication 
by $\lambda (1-\lambda)^{m-1}$.

\renewcommand{\theequation}{\thesection.\arabic{equation}}
\renewcommand{\therema}{\thesection.\arabic{rema}}
\setcounter{equation}{0}
\setcounter{rema}{0}

\section{\bf A Family Invariants of Unlabeled Posets}\label{S4}

In this section, we consider the invariants of unlabeled posets
derived from Algorithm \ref{Algo} by identifying unlabeled posets
with certain labeled posets.

First, let us identity unlabeled posets with naturally labeled posets. 
By Lemma \ref{L2.1}, it is easy to see that the restriction of
Algorithm \ref{Algo} on naturally labeled posets gives 
the following algorithm for order polynomials of unlabeled posets.

\begin{Algo}\label{N-Algo}
$(1)$ For the empty poset $\emptyset$, we set 
$\Psi (\emptyset) =a$.

$(2)$  For any non-empty poset $P$, 
we define  $\Psi (P)$ recursively by
\BQ
\Psi (P)=\Xi \sum_{\emptyset \neq I\in \cI (P)} \,   \Psi ({P\backslash S}).
\EQ
\end{Algo}

But on the other hand, if we identify unlabeled posets $P$
with labeled posets $(\Pw)$
with order-reversing labelings $\omega$, 
we get the following algorithm.

\begin{Algo}\label{UN-Algo}
For the empty poset $\emptyset$, we set 
$\Psi (\emptyset) =a$.

$(2)$  For any non-empty poset $P$, 
we define  $\Psi (P)$ recursively by
\BQ
\Psi (P)=\Xi \sum_{\emptyset \neq S\subset L(P)} \,  \Psi ({P\backslash S}).
\EQ
\end{Algo}

It is easy to see that, with $A=\bC[t]$, $a=1$ and 
$\Xi=\Delta$; the invariants defined by Algorithm \ref{N-Algo} and 
\ref{UN-Algo} are same as 
order polynomials $\Omega (P)$ and strict  
order polynomials $\bar \Omega (P)$
of unlabeled posets $P$, respectively.
Note that Algorithm \ref{UN-Algo} in general is much more 
efficient than Algorithm \ref{N-Algo}.
Order polynomials $\Omega (P)$ actually
can also be calculated  by 
Algorithm \ref{UN-Algo} as follows.

Let $\nabla$ be the $\bC$-linear operator defined 
\BQ
(\nabla f)(t)=f(t)-f(t-1).
\EQ
We also define the operator $\nabla^{-1}$ by setting  $(\nabla^{-1}f)(t)$ to 
be the unique polynomial $g(t)$ such that $\nabla f(t)=g(t)$ and $g(0)=0$.
From the definition of strict order polynomials, it is easy to see  
that 
\BQn
\bar \Omega (P; 1)=\begin{cases} 1 \quad \text{if $P$ is an anti-chain;}\\ 
0 \quad \text{otherwise.}\end{cases} 
\EQn
By the following well known identity (see \cite{St2}.)
\BQ\label{Recip}
\Omega (P; t)=(-1)^{|P|}\bar \Omega (P; -t),
\EQ
we have
\BQ\label{-1}
\Omega (P; -1)=\begin{cases} (-1)^{|P|} \quad \text{if $P$ is an anti-chain;}\\ 
0 \quad \text{otherwise.}\end{cases} 
\EQ

\begin{propo}
For any finite poset $P$ with $|P|\geq 2$, set $\wtilde \Omega (P; t)=(-1)^{|P|} \Omega (P; t)$. 
Then, we have
\BQ\label{nabla}
 \wtilde  \Omega (P; t)=-\nabla^{-1}
 \sum_{\emptyset \neq S\subset L(P)} \,  \wtilde \Omega (P\backslash S, t).
\EQ
\end{propo}

The proposition above can be proved by using Theorem \ref{s+t} and  Eq.\,(\ref{-1}).
It also can be proved by using Theorem \ref{Delta} and Eq.\,(\ref{Recip}).

Since $\wtilde \Omega (\emptyset, t)=1$, we have

\begin{propo}
Let $\Psi$ be the 
invariant 
defined by Algorithm \ref{UN-Algo} with
 $A=\bC [t]$, $a=1$ and $\Xi=-\nabla^{-1}$. 
Then, for any finite poset $P$, 
we have $\Psi_P=\wtilde \Omega (P; t)=
(-1)^{|P|} \Omega (P; t)$.
\end{propo}

In particular, we see that the order polynomials $\Omega (P; t)$ can also 
be recovered as an invariant of the type defined by Algorithm \ref{UN-Algo}.

{\small \sc Department of Mathematics, Washington University in St. Louis,
St. Louis, MO 63130.}

{\em E-mail}: shareshi@math.wustl.edu, wright@math.wustl.edu.

{\small \sc *Department of Mathematics, Illinois State University, Normal, IL 61790-4520.}

{\em *E-mail}: wzhao@ilstu.edu.


\begin{thebibliography}{FLM2}

\bibitem[B]{B} N. Biggs, {\it Algebraic Graph Theory}, 
Cambridge University Press, 1974.

\bibitem[G]{G} I. Gessel, {\it Multipartite $P$-partitions and
and Inner Products of Skew Schur Function}, Comtemp. Math. {\bf 34} (1984),
289-301.


\bibitem[MR]{MR} C. Malvenuto and C. Reutenauer, {\it Duality between 
Quasi-Symmertric  Functions and the Solomon Descent Algebra}, 
J. Algebra {\bf 177} (1995), 967-982.

\bibitem[St1]{St1}  Richard P. Stanley, 
{\it Ordered Structures and Partitions}, Ph.D Thesis, 
Harvard University, 1971. (See also
{\it Mem. Amer. Math, Soc.} {\bf 119} (1972).)

\bibitem[St2]{St2}  Richard P. Stanley, {\it Enumerative Combinatorics
I},
Cambridge University Press, 1997.

\bibitem[St3]{St3}  Richard P. Stanley, {\it Enumerative Combinatorics
II}, Cambridge University Press, 1999.

\bibitem[T]{T}  G. P. Thomas, {\it Frames, Yong Tableaux and Baxter Sequences},
Adv. in math. {\bf 26} (1997), 275-289.

\bibitem[WZ]{WZ} D. Wright and W. Zhao,
{\it D-log and formal flow for analytic isomorphisms of n-space}, 
{\it Trans. Amer. Math. Soc.}, {\bf 355}, No. 8, 3117-3141. 
Also see math.CV/0209274. 

\bibitem[Z]{Z} W. Zhao, {\it A Family of Invariants of Rooted Forests}.
math.CO/0211095. Appearing in {\it J. Pure Appl. Alg.} 


\end{thebibliography}
\end{document}